\documentclass[11pt,onecolumn]{article}
\usepackage{algorithm, algorithmic, setspace}
\usepackage{graphicx} % for pdf, bitmapped graphics files

\usepackage{amsmath} % assumes amsmath package installed
\usepackage{amssymb}  % assumes amsmath package installed
\usepackage{amsthm}
\usepackage{amsfonts}
\usepackage{dsfont}
\usepackage{color}
\usepackage{subfig}
\usepackage{booktabs}
\usepackage{subfig}
\usepackage{pgfplots}
\usepackage{tikz}
\usepackage[english]{babel}
\usepackage[latin1]{inputenc}

\newcommand{\xtrue}{\widetilde{x}}

\newcommand{\alf}{\mathcal{A}}

\newcommand{\R}{\mathds{R}}

\newtheorem{proposition}{Proposition}

\newtheorem{remark}{Remark}

\newtheorem{assumption}{Assumption}

\begin{document}
\title{Recovery of binary sparse signals from compressed linear measurements via polynomial optimization}

\author{Sophie M. Fosson and Mohammad Abuabiah }
\maketitle
\begin{abstract}
The recovery of signals with finite-valued components from few linear measurements is a problem with widespread applications and interesting mathematical characteristics. In the compressed sensing framework, tailored methods have been recently proposed to deal with the case of finite-valued sparse signals. In this work, we focus on  binary sparse signals and we propose a novel formulation, based on polynomial optimization. This approach is analyzed and compared to   the state-of-the-art binary compressed sensing methods.
\end{abstract}
\section{Introduction}
Solving underdetermined linear systems with finite-valued solutions has recently attracted a lot of attention in signal processing. This mathematical problem is encountered whenever a digital signal has to be recovered from compressed or incomplete linear measurements, which occurs in a number of applications, ranging from digital communications  and digital image processing, see, e.g., \cite{dra09,tia09,sto10,bio14,ahn16,hay17,hay18}, to localization, source separation, and spectrum sensing, see, e.g., \cite{pas15,lee16,kei17,liu18,fli18,fox19}.

Formally, the problem is defined as follows: given a finite alphabet $\alf$ (that is, a finite set of symbols), one aims to recover $\xtrue\in\alf^n$ from compressed measurements $y=A\xtrue$ (possibly corrupted by noise), where $A\in\R^{m,n}$, $m<n$, is a known  matrix. The solution is assumed to be unique within $\alf^n$. The problem is intrinsically combinatorial: if $c$ is the cardinality of $\alf$, the solution could be found by listing all the $c^n$ possible vectors, but clearly this is unfeasible for large $n$. 

If $\xtrue$ is sparse, i.e., it has many null components, sparsity-promoting procedures can be applied, in particular the ones developed in compressed sensing (CS, \cite{don06,fou13}). We recall that the CS paradigm states that (real-valued) sparse vectors can be recovered from few linear measurements, under suitable conditions on the matrix $A$ (for example, random matrices are fitting). Standard sparsity-promoting and CS techniques are conceived for real-valued signals, and do not envisage a possible prior information on the discrete nature of the signal. To fill this gap, in the last years, novel CS strategies have been proposed, which exploit the knowledge of $\alf$. These strategies are shown to significantly improve the performance with respect to algorithms unaware of the discrete structure, see, e.g., \cite{kei17,fli18,fox19} and the references therein. 

In the literature on the recovery of finite-valued sparse signals from linear compressed measurements, a lot of work is devoted to the binary case $\alf=\{0,1\}$. This line of research is also known as binary compressed sensing (BCS, \cite{nak12,shi15,fox18asi}). CS strategies tailored for BCS generally require less measurements than classical CS techniques. Different approaches are proposed, e.g., Bayesian models in \cite{tia09}, bi-partite graph models in \cite{nak12}, $\ell_0$ minimization in \cite{liu18}, $\ell_{\infty}$ minimization in \cite{ahn16},  $\ell_1$ minimization in \cite{sto10,pas15}, greedy methods in  \cite{fli18}, and non-convex models in \cite{fox18asi}. 

The interest on binary sparse signals is motivated by several applications. For example, in the framework of indoor localization, this setting is frequently considered \cite{bay15}: an area is subdivided into $n$ cells, and $m$ sensors are deployed to receive signals from targets on the area. In the training phase, a target transmits in turn from each cell, in order to build a dictionary $A\in\R^{m,n}$. Each entry $A_{i,j}$ represents the received signal strength measured by sensor $i$ when the target is in the cell $j$. In the running phase, the cells which contain transmitting targets are detected as solution of $y=Ax$, where $x\in\{0,1\}^n$. $A$ has then a mixing effect on signals transmitted by targets. Since usually the number of targets is much smaller than the number of cells, $x$ is sparse and CS can be applied. Beyond localization, the literature on the recovery of binary vectors from linear measurements is extensive in digital data transmission, jump systems, source separation, fault detection, see, e.g., \cite{tal96,fox13,fox15, beh17} and references therein. %In these works, a probabilistic approach is often proposed, which assumes the knoweldge of a prior distribution on the symbols.
Interestingly, we notice that the problem $Ax=y$ with $x\in\{0,1\}^n$  is also a variant of the popular 0-1 knapsack problem \cite{fay82}. Each equation of the system is a 0-1 knapsack problem whose goal is to fill a knapsack of capacity $y_i$ choosing among $n$ items with weights $A_{i,1},\dots, A_{i,n}$. %The 0-1 knapsack problem is known to be NP-hard, while pseudo-polynomial dynamic programming algorithms are available for the sub-case of positive integer weights.

In this work, we propose a novel polynomial optimization approach to recover binary sparse signals from compressed linear measurements. Specifically, we show that the problem can be formulated as a non-convex polynomial optimization problem (POP, \cite{las15}), whose solution can be estimated efficiently via semidefinite programming.
Our contribution can be summarized as follows. First, we propose a non-convex POP whose global minimum is exactly the desired solution, in the noise-free case. Second, we show that this POP can be solved via semidefinite programming; moreover, it has a sparse structure as described in  \cite{lasspa}, which reduces the computational complexity. Finally, we consider the noisy case, and we perform numerical tests to illustrate the recovery performance. In particular, we compare the proposed approach to the state-of-the-art BCS methods proposed in \cite{fli18,fox18asi}.

The paper is organized as follows. In Section \ref{sec:ps}, we introduce and analyze the proposed optimization problem. In Section \ref{sec:noise}, we extend the model to the noisy setting. In  Section \ref{sec:sims} is devoted to numerical simulations. Finally, we draw some conclusions.
\section{Problem statement and analysis}\label{sec:ps}
Let us consider the problem:
\begin{equation}\label{m1}
%\begin{split}
y=A x,~~~ x\in\{0,1\}^n,~~y\in\R^m,~A\in\R^{m,n},~~m< n
\end{equation}
where $A$ and $y$ are known. We assume that \eqref{m1} is well-posed.
\begin{assumption}\label{ass}
The solution $\xtrue\in\R^n$ of problem \eqref{m1} is unique. 
\end{assumption}
This assumption holds for many classes of $A$, for example, when the entries of $A$ are in general position (see \cite{tib13} for the definition), which is guaranteed when the entries of $A$ are generated at random according to any continuous distribution.

The key idea of this work is to tackle problem \eqref{m1} by reformulating it as follows:
\begin{equation}\label{pop}
\begin{split}
&\min_{x\in[0,1]^n} \sum_{i=1}^n\left(x_i-x_i^2\right)~~\text{ s. t. }~ y=A x.
\end{split}
\end{equation}
The optimization problem in \eqref{pop} is well-posed, as proven in the following proposition. 
\begin{proposition}\label{prop:1}
Under Assumption \ref{ass}, the unique solution of problem \eqref{pop} is $\xtrue$, for any $m\geq 1$.
\end{proposition}
\begin{proof}
The global minima of $f(x):=\sum_{i=1}^n\left(x_i-x_i^2\right)$ over $[0,1]^n$ are in $\{0,1\}^n$. In particular, $f(x)\geq 0$ for any $x\in[0,1]^n$, and $f(x)=0$ only when $x\in\{0,1\}$. This is consistent with the constraint $y=A x$, which has a unique binary solution $\xtrue$. Therefore, $\xtrue$ is the global minimum of problem \eqref{pop}. 
\end{proof}
Problem \eqref{pop} is non-convex, therefore its solution is not straightforward. However, since it is a POP, results proven in \cite{las01,par03,las15} can be exploited to the purpose. In a nutshell, these results state that the solution of a POP can be achieved by solving a sequence of relaxed semidefinite programming (SDP) problems, whose global minima converge to the global minimum of the POP. Thus, an estimate of the POP solution can be obtained by solving an SDP problem of sufficiently large relaxation order. On the other hand, the dimensions of the SDP problems increase with the relaxation order, therefore it is important to keep the order as low as possible. Interestingly, for Problem \eqref{pop}, we can prove that relaxation order 1 is sufficient to achieve the exact global minimum. 
\begin{proposition}\label{orderone}
The global minimum of the SDP relaxation of order 1 of Problem \eqref{pop} corresponds to the global minimum of Problem \eqref{pop}. The global minimizer of Problem \eqref{pop} can be extracted if the solution of the SDP relaxation has rank 1.
\end{proposition}
\begin{proof}
The first part follows from Theorem 4.2 in \cite{las01}, which, among others, states the following result. Let $p(x):\R^n\to\R$ be a polynomial and $K$ be the compact set described by $\{g_i(x)\geq 0, i=1,\dots,m\}$, where $g_i$'s are polynomials.  Let $p^{\star}=\min_{x\in K}p(x)$. If 
%\begin{equation}\label{SOScombine}
$p(x)-p^{\star}=q_0+\sum_{i=1}^m g_i(x) q_i$ 
%\end{equation}
where $q_j$, $j=0,\dots,m$, are non-negative constants, then the solution of the relaxed SDP of order $1$ is equal to $p^{\star}$.

For Problem \eqref{pop}, we know that $p^{\star}=0$, and we can replace $x_i\in [0,1]$ with $x_i-x_i^2\geq 0$, thus the condition of \cite[Theorem 4.2]{las01} holds with $q_0=0$, $q_i=1$ for $g_i(x)=x_i-x_i^2$, and $q_i=0$ for the constraints $Ax=b$.
Finally, the possibility of extracting the minimizer when the rank is 1 is a direct application of \cite[Theorem 6.6]{las15}.
\end{proof}
In conclusion, we can reformulate Problem \ref{pop} as a single SDP problem. More precisely, this SDP problem corresponds to the SDP relaxation for quadratic problems, which is written explicitly, e.g., in \cite[Section 4.5]{wak06}.
\begin{remark}
Problem \eqref{pop} and the related analysis can be extended to any binary alphabet $\{\alpha,\beta\}$ by using $(x-\alpha)(x-\beta)$ as cost functional. Similarly, extensions to any finite alphabet $\{\alpha_1,\dots,\alpha_c\}$, with $\alpha_1<\dots<\alpha_c$ are possible by using the cost functional $\sum_{i=1}^{c-1}[(x-\alpha_i)(x-\alpha_{i+1})]_+$, where $[\cdot]_+$ denotes the positive part. In this last case, the problem is semi-algebraic, and the theory in \cite{las15} still applies. Such extensions will be object of future work.
\end{remark}
%
%\vskip-2cm
\subsection{Complexity reduction via chordal sparsity}
Proposition \ref{orderone} recasts the non-convex problem \eqref{pop} into SDP, thus its solution can be computed via convex optimization. A drawback of this approach might be the dimension of the SDP problem, which increases with $n$. However, the \emph{chordal sparsity} property reduces the complexity, as described in \cite{lasspa,wak06,fan19}. Briefly, chordal sparsity is defined as follows: if the matrices involved in an SDP problem are sparse, the sparsity pattern can be represented by an undirected graph, see, e.g., \cite{fan19}. If this graph is chordal, i.e., each cycle composed by at least four nodes has a chord, then the SDP problem can be decomposed into smaller sub-problems (see Theorems 1 and 2 in \cite{fan19}), which reduces the complexity. 
 
 As illustrated in \cite{lasspa}, chordal sparsity  can be evaluated directly on the original POP as follows. Given a constrained POP with $n$ variables, such as \eqref{pop}, let us build $p$ sets $I_l\subset\{1,\dots,n\}$, $l=1,\dots,p$, with the following properties: their union is $\{1,\dots,n\}$; the variables $x_i$ contained in each constraint equation are concerned with a single $I_l$; the objective function is sum of monomials such that each monomial is concerned with a single $I_l$.  We then define the running intersection property as in \cite{lasspa}: for every $l=1,\dots,p-1$, $I_{l+1}\cap\left(\cup_{j=1}^l I_j\right)\subseteq I_s $ for some $s\leq l$. If the considered POP satisfies the running intersection property, then the graph associated with its SDP relaxation is chordal, and Theorems 1 and 2 in \cite{fan19} hold. In particular, the $I_l$'s represent the maximal cliques of such graph. 

Problem \eqref{pop} is sparse, as there are no mixed products between variables; this implies that the matrices in the SDP relaxation are null except for the main diagonal, the first row, and the first column, see \cite[Section 4.5]{wak06} for details. Moreover, Problem \eqref{pop} fulfills the running intersection property if we rewrite the constraints $Ax=y$ by adding suitable slack variables $z_{i,j}$, $i=1\dots,m$, $j=1\dots,n-1$. Specifically,
%
%{\small{
\begin{equation}\label{sparsepop}
\begin{split}
\min_{x\in[0,1]^n,z\in\R^{m,n-1}} \sum_{i=1}^n\left(x_i-x_i^2\right)&\text{ s. t., for each } i=1\dots,m, \\
y_i-A_{i,1} x_1+z_{i,1}&= 0\\
-A_{i,2}x_2-z_{i,1}+z_{i,2}&= 0\\
&\vdots\\
-A_{i,n}x_n-z_{i,n-1}&= 0.
\end{split}
\end{equation}
%}}
%We notice that by summing the constraints, the $z_{i,j}$'s are mutually canceled, and the original constraints in \eqref{pop} are retrieved. 
The following result holds (the proof is straightforward and omitted for brevity).
\begin{proposition} Problem \eqref{sparsepop} satisfies the running intersection property with the following sets: $I_{1}=\{x_1,z_{1,1},\dots,z_{m,1}\};$  $I_{l}=\{x_l,z_{1,l-1},z_{1,l},\dots, z_{m,l-1},z_{m,l}\}$ for $l=2,\dots,n-1$; $I_{n}=\{x_n,z_{1,n-1},\dots,z_{m,n-1}\}$.
\end{proposition}
According to \cite{lasspa,fan19}, the running intersection property allows us to recast our SDP problem of dimension $n^2$ into $p=n$ sub-problems of dimension $m^2$, which is favorable when $m\ll n$.

\subsection{Relation to prior literature}
We review the most important BCS approaches to highlight the differences with respect to the proposed Problem \eqref{pop}.
%In classical CS, the basis pursuit is very popular as convex relaxation of the problem, which reads as $\min_{x\in\R^n}\|x\|_1$ subject to $Ax=y$.
In \cite{sto10}, the authors propose to tackle BCS by solving a basis pursuit restricted to $[0,1]^n$, namely, 
$\min_{x\in[0,1]^n} f(x)=\sum_{i=1}^n x_i$ such that  $y=A x$, $A\in\R^{m,n}$, $m<n$.
This model, later revisited in \cite{pas15,kei17}, is shown to require less measurements than classical basis pursuit, and, in particular, when $m> n/2$, it provides the exact solution with high probability for specific random $A$, see  \cite{sto10,pas15}. A similar result is obtained by using linear programming in \cite{man11}. In \cite{hay17,hay18}, the authors propose a sum-of-absolute-value optimization by assuming to know the probability of $x_i=0$. % which reads as follows: $\min p\|x\|_1+(1-p)\|x-(1,\dots,1)^T\|_2$ subject to $Ax=y$.  
In \cite{fox18asi}, a non-convex model is proposed: $\min_{x\in[0,1]^n}\frac{1}{2}\left\|y-A x \right\|_2^2+\lambda\sum_{i=1}^n\left(x_i-\frac{x_i^2}{2}\right)$, $\lambda>0$, which can be interpreted as a Lasso \cite{tib96} with concave (non-decreasing) penalization instead of $\ell_1$ penalization to compel sparsity. The global minimum is exactly the desired signal, under weak conditions, and $\ell_1$ iterative reweighting algorithms are  shown to achieve it, if $m$ is sufficiently large. In \cite{fli18}, the PROMP algorithm, based on orthogonal matching pursuit, is proposed, which requires a larger $m$ than \cite{fox18asi}, while it is faster.
\section{Extension to the noisy case}\label{sec:noise}
Problem \eqref{pop} can be extended to the case of noisy measurements, namely $y=A\xtrue+\zeta$, $\zeta\in\R^n$ being an unknown disturbance. If we assume that $\|\zeta\|_{\infty}\leq \eta$, where $\eta>0$ is known, Problem \eqref{pop} can be reformulated as follows:
\begin{equation}\label{popnoise}
\begin{split}
&\min_{x\in[0,1]^n} \sum_{i=1}^n\left(x_i-x_i^2\right)\text{ s. t. } \|Ax-y\|_{\infty}\leq \eta.
\end{split}
\end{equation}
As $\|Ax-y\|_{\infty}\leq \eta$  corresponds to $|A_i x-y_i|\leq \eta $ for each $i=1,\dots,n$, we can write the sparse formulation as follows:
%
%{\small{
\begin{equation}\label{sparsepopnoise}
\begin{split}
\min_{x\in[0,1]^n,z,w\in\R^{m,n-1}} \sum_{i=1}^n\left(x_i-x_i^2\right)&\text{ s. t. for each } i=1\dots,m, \\
\eta+y_i-A_{i,1} x_1+z_{i,1}\geq 0,~~~&\eta-y_i+A_{i,1} x_1+w_{i,1}\geq 0\\
-A_{i,2}x_2-z_{i,1}+z_{i,2}\geq 0,~~~&A_{i,2}x_2-w_{i,1}+w_{i,2}\geq 0\\
&\vdots\\
-A_{i,n}x_n-z_{i,n-1}\geq 0,~~~&A_{i,n}x_n-w_{i,n-1}\geq 0
\end{split}
\end{equation}
%}}
where $z_{i,j}$ and $w_{i,j}$, $i=1\dots,m$, $j=1\dots,n-1$, are slack variables. The number of slack variables is doubled with respect to the noise-free case, which increases the whole dimension of the problem. However, since the constraints are all linear, the increase of complexity is not dramatic. Moreover, if $\eta$ is sufficiently small, we can still expect a unique solution under the same conditions of the noise-free case.
\section{Numerical results}\label{sec:sims}
In this section, we propose numerical simulations that support the effectiveness of the proposed approach.
We consider the following setting: the unknown binary signal $\xtrue\in\{0,1\}^n$ has $n=100$ components, among which  $k\in[10,50]$ components are equal to 1, with  sparsity ratio $\frac{k}{n}\in\left[\frac{1}{10}, \frac{1}{2}\right]$. We test $m\in[20,50]$. The support of the signal is generated uniformly at random, and the entries of $A$ are Gaussian $\mathcal{N}(0,\frac{1}{m})$. To recover $\xtrue$, we solve the POP's \eqref{pop} and \eqref{popnoise} by solving the corresponding SDP relaxations with the solver CDCS, see \cite{fan19,cdcs}, on Matlab R2018a. CDCS solves SDP problems by implementing the alternating direction method of multipliers (ADMM, \cite{boy10}), which is faster than, e.g., interior point methods. Moreover, CDCS exploits the chordal sparsity of SDP to efficiently reduce the dimensionality of the problem, which perfectly matches with problems \eqref{sparsepop} and \eqref{sparsepopnoise}. 

Our approach is compared to the state-of-the-art BCS algorithms RWR \cite{fox18asi} and PROMP \cite{fli18}. RWR is based on repeating a local minimization iterative algorithm starting from different random initial points, and is shown to outperform basis pursuit based methods, see \cite[Section IV]{fox18asi}.
The following recovery accuracy metrics are considered (we denote by $\widehat{x}$ the estimate of $\xtrue$): exact recovery rate, i.e., how many times the whole signal is exactly recovered; false positive rate, i.e., the number of occurrences $\widehat{x}_i \neq 0$ when $\xtrue_i=0$ divided by the true number of zeros $n-k$; false negative rate, i.e., the number of occurrences $\widehat{x}_i = 0$ when $\xtrue_i \neq 0$ divided by the true number of ones $k$.
The run time is illustrated as well. The results are averaged over 100 runs. Estimations are quantized to $\{0,1\}$ (by imposing a threshold at $\frac{1}{2}$) when algorithms do not provide binary solutions.
\subsection{Noise-free experiment}
\begin{figure}[ht]
\centering
\includegraphics[width=1\columnwidth]{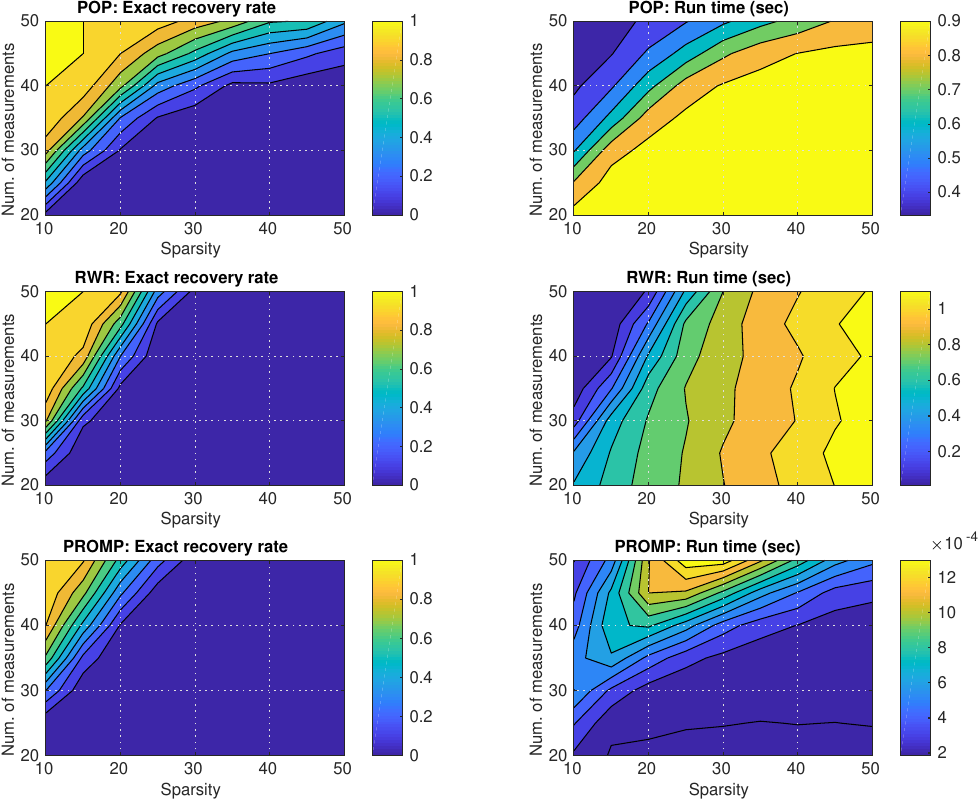}\caption{Noise-free case: exact recovery rate and run time (in seconds); $n=100$, $k\in[10,50]$, $m\in[20,50]$.}\label{fig1}
\end{figure}
\begin{figure}[ht]
\centering
\includegraphics[width=1\columnwidth]{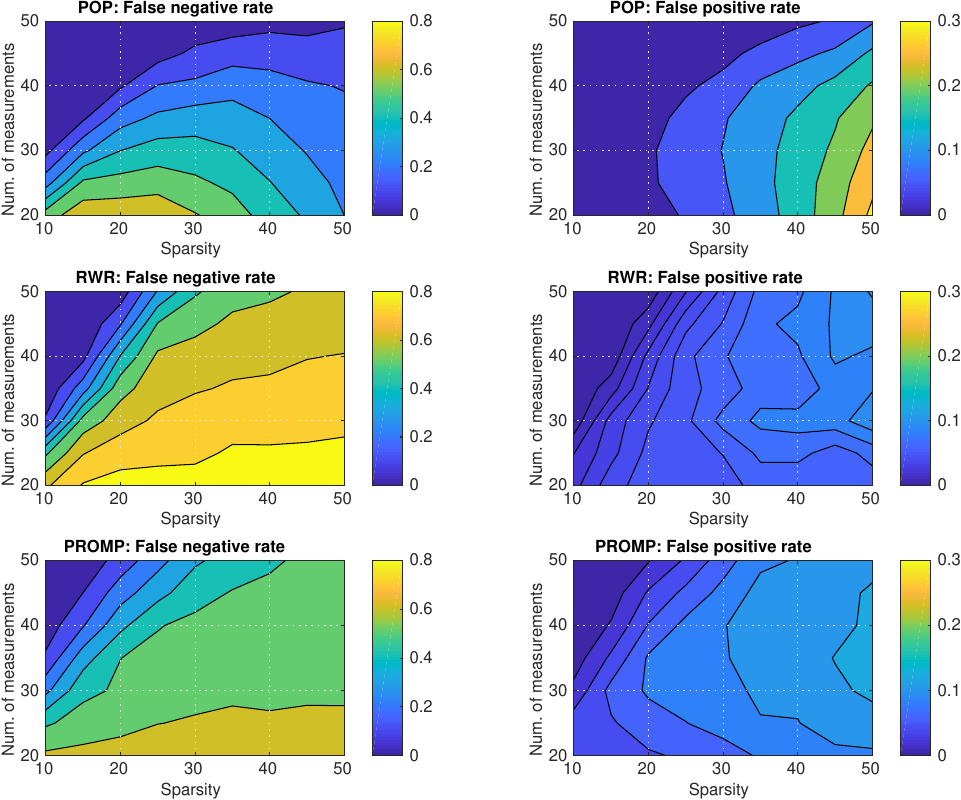}\caption{Noise-free case: false negatives and false positives.}\label{fig2}
\end{figure}

In Figures \ref{fig1} and \ref{fig2}, we show the results in the  noise-free case. In Figure \ref{fig1}, we can see that the proposed POP approach provides the best recovery accuracy, while PROMP is the least accurate. In particular, in contrast to RWR and PROMP, for POP $m=45$ is always sufficient to recover exactly a signal with $k\leq 15$ (the light yellow denotes $100\%$ of success). We remark that the ratio $\frac{m}{k}=3$ is favorable with respect to CS standards \cite{fou13}. When $k\geq 30$, RWR and PROMP always fail ($0\%$ of success), while POP is still successful in more than $50\%$ of runs with $m=50$. POP and RWR have similar run time, of order $10^{-1}$ seconds. PROMP is faster ($10^{-4}$ seconds); nevertheless, as already said, its performance accuracy is scarce.

From Proposition \ref{prop:1}, one might expect to find the exact solution even when $m=1$, using the POP approach. However, the corresponding SDP problem might have several global minima. To obtain the right minimizer the SDP solution must have rank 1. A way to induce low rank could the penalization of the nuclear norm, see \cite{rec10}; this strategy will be investigated in the future. In this work, we observe that, in practice, a sufficient small $k$ and sufficient large $m$ provide the exact solution.

In Figure \ref{fig2}, we further examine the estimation accuracy by considering the correct recovery of single symbols, in terms of false negative and false positive rate. POP has a more favorable false negative rate with respect to RWR and PROMP, as it does not not induce sparsity. On the other hand, the false positive rate is similar for the three methods, except for large $k$. 
\subsection{Noisy experiment}
We finally revisit the previous experiment by adding measurement noise. We set $\eta=5\times 10^{-2}$, which corresponds to a signal-to-noise ratio $\frac{k}{m\eta^2}$ between $25$ and $40$ dB.
\begin{figure}[ht]
\centering
\includegraphics[width=1\columnwidth]{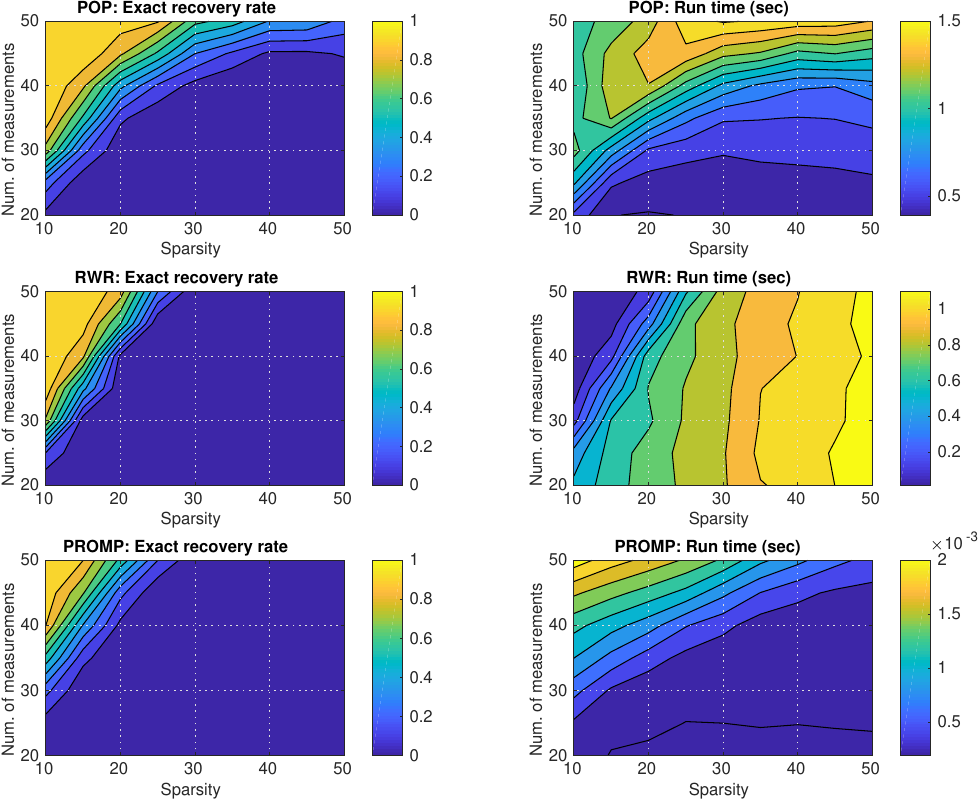}\caption{Noisy case: exact recovery rate and run time (in seconds); $n=100$, $k\in[10,50]$, $m\in[20,50]$.}\label{fig1n}
\end{figure}
\begin{figure}[ht]
\centering
\includegraphics[width=1\columnwidth]{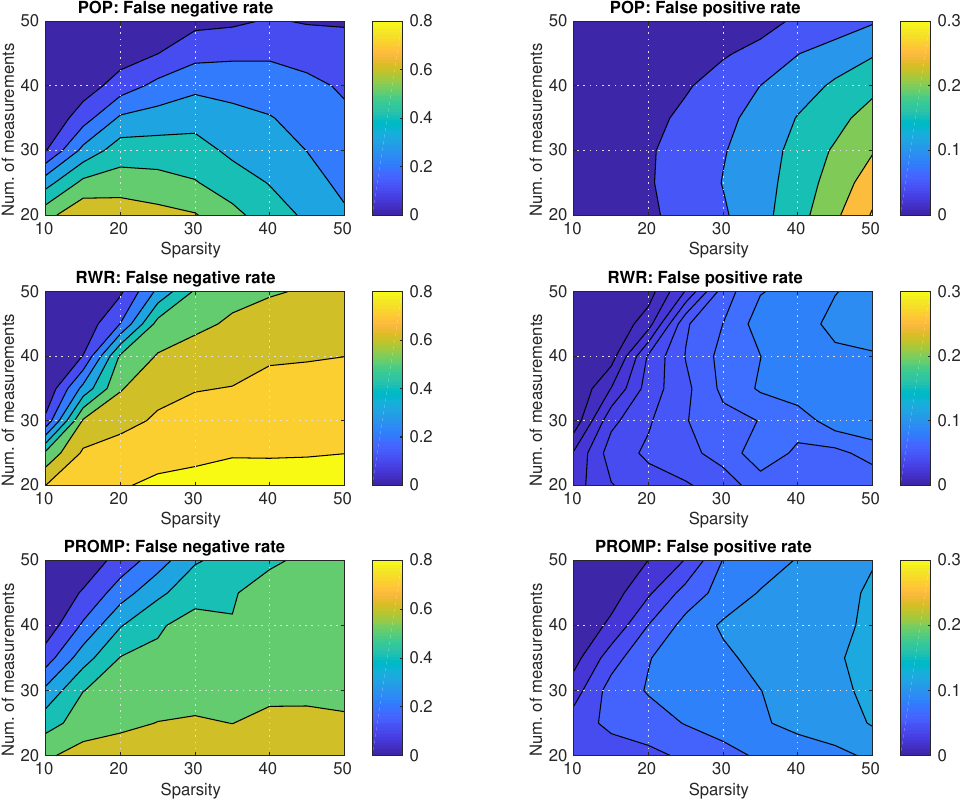}\caption{Noisy case: false negatives and false positives.}\label{fig2n}
\end{figure}
In Figures \ref{fig1n} and \ref{fig2n}, we can see that, as expected, the general performance is slightly worse when noise is present. However, the mutual behaviors between POP, RWR, and PROMP are similar to those observed in the noise-free case.
\vspace{-0.05cm}
\section{Conclusion}
In this paper, we develop a non-convex quadratic approach to recover binary signals from few linear measurements. The effectiveness of the proposed method leverages on recent results on polynomial optimization, and has better recovery accuracy than binary compressed sensing algorithms, in particular when signals are not very sparse. Future work will envisage the extension of the proposed approach to larger, non-binary, alphabets, and to non-sparse problems, e.g., combinatorial problems, like the 0-1 knapsack problem.
\bibliographystyle{plain}
\bibliography{refs}

\begin{thebibliography}{10}

\bibitem{ahn16}
J.~H. Ahn.
\newblock Compressive sensing and recovery for binary images.
\newblock {\em IEEE Trans. Image Process.}, 25(10):4796--4802, 2016.

\bibitem{pas15}
A.~A\"{\i}ssa-El-Bey, D.~Pastor, S.~M.~A. Sba\"{\i}, and Y.~Fadlallah.
\newblock Sparsity-based recovery of finite alphabet solutions to
  underdetermined linear systems.
\newblock {\em IEEE Trans. Inf. Theory}, 61(4):2008--2018, 2015.

\bibitem{bay15}
A.~Bay, D.~Carrera, S.~M. Fosson, P.~Fragneto, M.~Grella, C.~Ravazzi, and
  E.~Magli.
\newblock Block-sparsity-based localization in wireless sensor networks.
\newblock {\em EURASIP J. Wirel. Commun. Netw.}, 2015(182):1--15, 2015.

\bibitem{beh17}
M.~Behr and A.~Munk.
\newblock Identifiability for blind source separation of multiple finite
  alphabet linear mixtures.
\newblock {\em IEEE Trans. Inf. Theory}, 63(9):5506--5517, 2017.

\bibitem{bio14}
V.~Bioglio, G.~Coluccia, and E.~Magli.
\newblock Sparse image recovery using compressed sensing over finite alphabets.
\newblock In {\em IEEE Int. Conf. Image Process. (ICIP)}, pages 1287--1291,
  2014.

\bibitem{boy10}
S.~Boyd, N.~Parikh, E.~Chu, B.~Peleato, and J.~Eckstein.
\newblock Distributed optimization and statistical learning via the alternating
  direction method of multipliers.
\newblock {\em Found. Trends Mach. Learn.}, 3(1):1 -- 122, 2010.

\bibitem{don06}
D.~L. Donoho.
\newblock Compressed sensing.
\newblock {\em IEEE Trans. Inf. Theory}, 52(4):1289--1306, 2006.

\bibitem{dra09}
S.~C. Draper and S.~Malekpour.
\newblock Compressed sensing over finite fields.
\newblock In {\em Proc. IEEE Int. Symp. Inf. Theory (ISIT)}, pages 669--673,
  2009.

\bibitem{fox15}
F.~Fagnani and S.~M. Fosson.
\newblock Analysis of reduced-search bcjr algorithms for input estimation in a
  jump linear system.
\newblock {\em Signal Process.}, 108:341 -- 350, 2015.

\bibitem{fay82}
D.~Fayard and G~Plateau.
\newblock An algorithm for the solution of the 0-1 knapsack problem.
\newblock {\em Computing}, (28), 1982.

\bibitem{fli18}
A.~Flinth and G.~Kutyniok.
\newblock {PROMP}: A sparse recovery approach to lattice-valued signals.
\newblock {\em Appl. Comput. Harmon. Anal.}, 45(3):668--708, 2018.

\bibitem{fox13}
S.~M. Fosson.
\newblock Binary input reconstruction for linear systems: A performance
  analysis.
\newblock {\em Nonlin. Anal. Hybr. Syst.}, 7(1):54 -- 67, 2013.

\bibitem{fox18asi}
S.~M. Fosson.
\newblock Non-convex approach to binary compressed sensing.
\newblock In {\em Asilomar Conf. Signals Syst. Comput.}, 2018.

\bibitem{fox19}
S.~M. Fosson.
\newblock Non-convex {L}asso-kind approach to compressed sensing for
  finite-valued signals.
\newblock {\em (under review)}, 2019.

\bibitem{fou13}
S.~Foucart and H.~Rauhut.
\newblock {\em A Mathematical Introduction to Compressive Sensing}.
\newblock Springer, New York, 2013.

\bibitem{hay17}
R.~Hayakawa and K.~Hayashi.
\newblock Convex optimization-based signal detection for massive overloaded
  {MIMO} systems.
\newblock {\em IEEE Trans. Wireless Commun.}, 16(11):7080--7091, 2017.

\bibitem{hay18}
R.~Hayakawa and K.~Hayashi.
\newblock Discreteness-aware approximate message passing for discrete-valued
  vector reconstruction.
\newblock {\em IEEE Trans. Signal Process.}, 66(24):6443--6457, 2018.

\bibitem{kei17}
S.~Keiper, G.~Kutyniok, D.~Gwan Lee, and G.~E. Pfander.
\newblock Compressed sensing for finite-valued signals.
\newblock {\em Lin. Alg. Appl.}, 532(Supplement C):570--613, 2017.

\bibitem{las01}
J.~Lasserre.
\newblock Global optimization with polynomials and the problem of moments.
\newblock {\em SIAM J. Optim.}, 11(3):796--817, 2001.

\bibitem{lasspa}
J.~B. Lasserre.
\newblock Convergent {SDP}-relaxations in polynomial optimization with
  sparsity.
\newblock {\em SIAM J. Optim.}, (17):822--843, 2006.

\bibitem{las15}
J.~B. Lasserre.
\newblock {\em An Introduction to Polynomial and Semi-Algebraic Optimization}.
\newblock Cambridge Texts in Applied Mathematics. Cambridge University Press,
  2015.

\bibitem{lee16}
N.~Lee.
\newblock {MAP} support detection for greedy sparse signal recovery algorithms
  in compressive sensing.
\newblock {\em IEEE Trans. Signal Process.}, 64(19):4987--4999, 2016.

\bibitem{liu18}
T.~Liu and D.~G. Lee.
\newblock Fast binary compressive sensing via $\ell_0$ gradient descent.
\newblock 2018.

\bibitem{man11}
O.L. Mangasarian and B.~Recht.
\newblock Probability of unique integer solution to a system of linear
  equations.
\newblock {\em Eur. J. of Oper. Res.}, 214(1):27--30, 2011.

\bibitem{nak12}
U.~Nakarmi and N.~Rahnavard.
\newblock {BCS}: Compressive sensing for binary sparse signals.
\newblock In {\em IEEE Milit. Commun. Conf. (MILCOM)}, pages 1--5, 2012.

\bibitem{par03}
P.~A. Parrilo.
\newblock Semidefinite programming relaxations for semialgebraic problems.
\newblock {\em Math. Program.}, 96(2):293--320, 2003.

\bibitem{rec10}
B.~Recht, M.~Fazel, and P.~Parrilo.
\newblock Guaranteed minimum-rank solutions of linear matrix equations via
  nuclear norm minimization.
\newblock {\em SIAM Review}, 52(3):471--501, 2010.

\bibitem{shi15}
M.~Shirvanimoghaddam, Y.~Li, B.~Vucetic, J.~Yuan, and P.~Zhang.
\newblock Binary compressive sensing via analog fountain coding.
\newblock {\em IEEE Trans. Signal Process.}, 63(24):6540--6552, 2015.

\bibitem{sto10}
M.~Stojnic.
\newblock Recovery thresholds for $\ell_1$ optimization in binary compressed
  sensing.
\newblock In {\em Proc. IEEE Int. Symp. Inf. Theory (ISIT)}, pages 1593--1597,
  2010.

\bibitem{tal96}
S.~Talwar, M.~Viberg, and A.~Paulraj.
\newblock Blind separation of synchronous co-channel digital signals using an
  antenna array. i. algorithms.
\newblock {\em IEEE Trans. Signal Process.}, 44(5):1184--1197, 1996.

\bibitem{tia09}
Z.~Tian, G.~Leus, and V.~Lottici.
\newblock Detection of sparse signals under finite-alphabet constraints.
\newblock In {\em IEEE Int. Conf. Acoust. Speech Signal Process. (ICASSP)},
  pages 2349--2352, 2009.

\bibitem{tib96}
R.~Tibshirani.
\newblock Regression shrinkage and selection via the {L}asso.
\newblock {\em J. Roy. Stat. Soc. Series B}, 58:267--288, 1996.

\bibitem{tib13}
Ryan~J. Tibshirani.
\newblock {The {L}asso problem and uniqueness}.
\newblock {\em Electronic Journal of Statistics}, 7:1456--1490, 2013.

\bibitem{wak06}
H.~Waki, S.~Kim, M.~Kojima, and M.~Muramatsu.
\newblock Sums of squares and semidefinite program relaxations for polynomial
  optimization problems with structured sparsity.
\newblock {\em SIAM J. Optim.}, 17(1):218--242, 2006.

\bibitem{cdcs}
Yang Zheng, Giovanni Fantuzzi, Antonis Papachristodoulou, Paul Goulart, and
  Andrew Wynn.
\newblock {CDCS}: Cone decomposition conic solver, version 1.1.
%\newblock \url{https://github.com/giofantuzzi/CDCS}, 2016.

\bibitem{fan19}
Yang Zheng, Giovanni Fantuzzi, Antonis Papachristodoulou, Paul Goulart, and
  Andrew Wynn.
\newblock Chordal decomposition in operator-splitting methods for sparse
  semidefinite programs.
\newblock {\em Math. Progamm.}, 2019.

\end{thebibliography}
\end{document}